\documentclass[10pt]{amsart}
\usepackage{latexsym}
\usepackage{amsmath}
\usepackage{amssymb}
\usepackage{mathrsfs}
\begin{document}

\title{
LINEAR PROGRAMMING AND KANTOROVICH SPACES}

\author{S.~S. Kutateladze}
\address[]{
Sobolev Institute of Mathematics\newline
\indent Novosibirsk 630090\newline
\indent Russia}
\email{
sskut@member.academ.org
}
\begin{abstract}

This is a brief overview of the life  of Leonid Kantorovich
(1912--1986) and his contribution to the fields of linear programming and ordered vector spaces.
\end{abstract}
\maketitle

{\baselineskip=.975\baselineskip
January 19, 2007 is the date of the 95th anniversary of the birth of
Leonid  Kantorovich.   Among the Russian scientists,
Leonid Kantorovich and Vasily Leontiev were the only persons awarded with
a Nobel prize in economics. Kantorovich occupies a rather special place
in the history of the world science, ranking among the giants
that synthesized the exact and humanitarian modes of thought in
their creative activities.
Kantorovich was elected a corresponding member of
the Academy of Sciences of the USSR in economics and  a full member
in mathematics. There are no persons with similar biographic signposts in
the national academy of this country.
Kantorovich was among the scientists that comprised the first staff of the
Siberian Division of the Academy of Sciences of the USSR and spent a decade in
Akademgorodok near Novosibirsk.

\section{The World  Line of Kantorovich}

Kantorovich was born in the~family of a~venereologist at St.~Petersburg on January 19, 1912 (January 6, according to the
old Russian style).
It is curious that many reference books give another date
(which is three days before).
Kantorovich  kept explaining with a~smile that
he remembers himself from January 19, 1912.
The boy's talent was revealed   very early.
In 1926, just  at the age of 14,  he entered St. Petersburg
(then Leningrad) State University (SPSU).  Soon he started participating
in a~circle of G.~M.~Fikhtengolts
for students
and in a~seminar on  descriptive function theory.
It is natural that the early academic years formed his first environment:
D.~K.~Faddeev, I.~P.~Natanson, S.~L.~Sobolev, S.~G.~Mikhlin, and a~few others
with
whom Kantorovich was friendly during all his life also participated in
Fikhtengolts's circle. The old cronies
called him ``L\"enechka'' ever since these days.

After graduation from SPSU in 1930, Kantorovich started
teaching, combining it with intensive scientific research. Already in
1932  he became a~full professor at  the Leningrad Institute of
Civil Engineering and an assistant professor at SPSU.
From 1934 Kantorovich
was a~full professor at his {\it alma mater}.

The main achievements in mathematics belong to the
``Leningrad'' period of Kantorovich's life. In the
1930s he published more papers in pure mathematics
whereas his 1940s are devoted to computational mathematics
in which he was soon appreciated as a~leader in this country.


The letter of Academician  N.~N. Luzin
\cite{Reshetnyak}, written on April 29, 1934,
was  found in the personal archive of Kantorovich a few years ago
during preparation of his selected works for publication \cite{Kantorovich}.

This letter demonstrates the attitude of Luzin, one of the most eminent
and influential mathematicians of that time, to the brilliance of
the young prodigy. Luzin was one the nest mathematicians of that epoch and
the founder of the famous ``Lusitania'' school of Muscovites.
He remarked in his letter:

\smallskip
\begin{itemize}
\item[]{\small\sl\indent
{\dots}you must know my attitude to you.  I do not know you as a man
$\underline{\text{completely}}$
but I guess a warm and admirable personality.

However, one thing I know for certain: the range of your mental powers
which, so far as I  accustomed myself to guess people, open up
$\underline{\text{limitless possibilities}}$  in science.   I will not utter
the appropriate word---what for?
Talent---this would belittle you. You are entitled to get more\dots.}
\end{itemize}
\smallskip

\noindent
From the end of the 1930s the research of Kantorovich acquires new traits as
made an outstanding breakthrough in economics.
His booklet
{\it Mathematical Methods in the Organization and Planning of Production\/}
that appeared in 1939 is a material evidence
of the birth of linear programming. The economic  works of
Kantorovich were hardly visible at the surface of the scientific
information flow in the 1940s. However, the problems of economics
prevailed in his creative studies. During the Second World War  he
completed the first version of his book {\it The Best Use of Economic
Resources} which led to the Nobel Prize awarded to him and
Tjalling C.~Koopmans in~1975.

In 1957 Kantorovich accepted the invitation to join
the newly founded  Siberian Division of the Academy of Sciences of
the USSR.
He agreed and soon became a corresponding member
of the Department of Economics in the first elections to the
Siberian Division.
Since then  his major publications were
devoted to economics, with the exception of
the celebrated course of
functional analysis \cite{Akilov},
``Kantorovich and Akilov'' in the students' jargon.

I cannot help but mention one brilliant twist of mind of
Kantorovich and his students in suggesting a scientific approach to
taxicab metered rates.
The people of the elder generation remember that in the 1960s the taxicab meter
rates in the USSR were modernized radically: there appeared a price for
taking a taxicab which was combining with a less  per kilometer cost.
This led immediately to raising efficiency of
taxi parks as well as profitability of short taxicab drives.
This economical measure was a result of a mathematical modeling
of taxi park efficiency which was accomplished
by Kantorovich with a group of young mathematicians and
published in the most prestigious mathematical journal
``Russian Mathematical Surveys.''

The 1960s became the decade of his recognition.
In 1964 he was elected a full member of the Department of Mathematics
of the Academy of Sciences of the USSR,
and in 1965 he was awarded the Lenin Prize. In these years he
vigorously propounded and maintained his views of interplay between
mathematics and economics  and exerted great efforts to instill the ideas and methods
of  modern science into the top economic management of the Soviet Union,
which was almost in vain.

At the beginning of the 1970s Kantorovich left Novosibirsk for
Moscow where he was still engaged in economic analysis, not
ceasing his efforts to influence  the everyday economic practice and
decision making  in the national economy. These years witnessed
a~considerable mathematical Renaissance of Kantorovich. Although
he never resumed  proving theorems, his impact on the mathematical
life of this country increased sharply due to the sweeping changes
in the Moscow academic life on the eve of Gorbi's ``perestro\u\i{}ka.''
Cancer terminated his path in science on April~7, 1986.
He was buried at Novodevichy Cemetery in Moscow.

\section{Contribution to Science}

The scientific legacy of Kantorovich is immense \cite{Kantorovich}. His research in
the areas of functional analysis, computational mathematics,
optimization, and descriptive set theory has had a dramatic impact
on the foundation and progress of these disciplines. Kantorovich
deserved his status of one of the father founders of the modern
economic-mathematical methods. Linear programming, his most
popular and celebrated discovery, has changed the image of economics.

Kantorovich wrote more than 300 articles. When we discussed
with him the first edition of an annotated bibliography of his
publications in the early 1980s, he  suggested to combine them in
the nine sections: descriptive function theory and set theory,
constructive function theory,
approximate methods of analysis,
functional analysis,
functional analysis and applied mathematics,
linear programming,
hardware and software,
optimal planning and optimal prices,
and the economic problems of a~planned economy.

The impressive diversity of these areas of research
rests upon not only  the traits of Kantorovich but also
his methodological views.
He always emphasized the innate
integrity of his scientific research as well as mutual penetration and
synthesis of the methods and techniques he used in solving the
most diverse theoretic and applied problems of mathematics and
economics. I leave a thorough analysis of the methodology of
Kantorovich's contribution a challenge to professional
scientometricians. It deserves mentioning right away only that the abstract
ideas of Kantorovich in the theory of Dedekind complete vector
lattices, now called {\it Kantorovich spaces\/} or $K$-{\it
spaces},\footnote{Kantorovich wrote about ``my spaces'' in his personal memos.}
turn out to be closely connected with the art of
linear programming and the approximate methods of analysis.

It is impossible to shed light
on everything invented by Kantorovich in a short article.
He told me once that his main mathematical achievement
is the development of the theory of $K$-spaces
within functional analysis while remarking that his most useful deed
is linear programming.
These beautiful pearls of the scientific legacy of Kantorovich
deserve a special overview.

\section{Linear Programming}

Linear programming belongs to the curricula
of  hundreds of thousands of students throughout the world.
This term signifies a colossal scientific area
that addresses linear models of optimization.
 In other words, the science and art  of linear programming  deal with
 a theoretical and numerical analysis as well as solution of
 the problems in which we seek for an optimal value, i.~e.
 a maximum or minimum, of some system of indicators in any process
 whose states and/or behavior are determined by
 simultaneous linear inequalities.
The term ``linear programming'' was suggested in 1951
by Koopmans who contributed greatly
to proliferation of the idea of linear programming
and defending the  priority of Kantorovich in this
area.

In the USA, linear programming started only in 1947  by
George B. Dantzig who convincingly described the history of linear programming
in his classical
book \cite[p.~22--23]{Dantzig}
as follows:

\smallskip
\begin{itemize}
\item[]{\small\sl\indent
The  Russian  mathematician L. V.  Kantorovich  has  for  a
number  of  years  been  interested in  the  application  of
mathematics   to  programming  problems.  He  published   an
extensive monograph in 1939 entitled {\it Mathematical Methods in
the Organization and Planning of Production}\dots.

Kantorovich  should  be credited with being  the  first  to
recognize that certain important broad classes of production
problems had well-defined mathematical structures which,  he
believed,  were  amenable to practical numerical  evaluation
and could be numerically solved.

 In  the  first  part of his work Kantorovich  is  concerned
with  what  we  now call the weighted two-index distribution
problems. These were generalized first to include  a  single
linear  side  condition,  then  a  class  of  problems  with
processes     having     several    simultaneous     outputs
(mathematically the latter is equivalent to a general linear
program). He outlined a solution approach based on having on
hand  an  initial feasible solution to the  dual.  (For  the
particular problems studied, the latter did not present  any
difficulty.)  Although the dual variables  were  not  called
``prices,''  the general idea is that the assigned  values  of
these  ``resolving multipliers'' for resources in short supply
can  be  increased  to a point where it  pays  to  shift  to
resources that are in surplus. Kantorovich showed on  simple
examples  how  to make the shifts to surplus  resources.  In
general,  however, how to shift turns out  to  be  a  linear
program  in  itself  for which no computational  method  was
given.  The  report  contains an outstanding  collection  of
potential applications\dots.

If  Kantorovich's earlier efforts had been  appreciated  at
the  time  they  were first presented, it is  possible  that
linear  programming would be more advanced  today.  However,
his  early work in this field remained unknown both  in  the
Soviet  Union  and  elsewhere for nearly two  decades  while
linear  programming became a highly developed art.}
\end{itemize}
\smallskip

\noindent
It is wort observing that to an optimal plan of every
linear program there corresponds some optimal prices
or ``objectively determined estimators.'' Kantorovich
invented this  bulky by tactical reasons in order to
enhance the ``criticism endurability'' of  the concept.
The interdependence of optimal solutions and optimal prices is
the crux of the economical discovery of Kantorovich.

\section{Kantorovich Spaces}

In the mid-1930s the research of Kantorovich  contributed to the foundation of
a new fruitful area of functional analysis, the theory of ordered vector
spaces. Kantorovich introduced and elaborated the class of the so-called
Dedekind complete vector lattices in which each bounded subset has
a~least upper bound. These spaces are often referred to as
$K$-spaces or Kantorovich spaces. Kantorovich
gave versatile applications of his theory to many problems of
functional analysis, function theory, and the theory of
operator equations.

Kantorovich incessantly emphasized the immanent connection between
the theory of $K$-spaces, the theory of inequalities, and economical topics.
The further research of various authors demonstrates that the idea
of linear programming are inseparable from the theory of $K$-spaces in the
following rigorous sense:
The validity of each of
the universally accepted formulations of the duality
principle with prices in some algebraic structure necessitates  that
this  structure is a~$K$-space.

The  development of Boolean valued model of set theory
began in the 1960s in connection with the Paul Cohen's final solution
of the continuum problem. This problem was raised by Davis Hilbert
as the firsts in the list of his epochal report to the
International  Congress of Mathematicians in 1900.
The progress of the recent Boolean valued analysis
\cite{IBA} has demonstrated the fundamental significance
of  universally complete Kantorovich spaces.
Unexpectedly, each of these spaces turns out
to be one of the legitimate models of the real line and so
plays the same key role in mathematics.
It is curious that the development of new logical tools
led in the 1980s to reinvention of  $K$-spaces in the USA
under the title ``Boolean linear spaces.''

Magically prophetic happens to be the claim of Kantorovich
that the elements of a~$K$-space are generalized numbers.
The heuristic transfer principle
of Kantorovich  found a~brilliant justification
in the framework of modern mathematical logic.
Guaranteeing a profusion of  unbelievable
models of the real axis, the spaces of Kantorovich
will stay in the treasure-trove of the world science for ever.

\section{Synthesis of Cultures}

Kantorovich is rightfully ranked among the
founders of the economical-mathema\-tic\-al area.
Linear programming, his best invention, changed the appearance
of economics.
Kantorovich was in full possession not only of the gift of a mathematical
genius but also stamina of an intellectual champion in
struggling for  new economical theories.

The ideas and methods of linear programming
gave rise to deep interdisciplinary research,
trespassed the frontiers of economics, and won
appreciation in the various spheres of human activities.
It is difficult to distinguish another scholar in the history
of the twentieth century science who contributed as much as him
to interpenetration of mathematics and economics, the science with
the antipodal standards of scientific thought.
Israel Gelfand pointed out that
he can list only John von Neumann and Andrei Kolmogorov
alongside Leonid Kantorovich among his contemporaries
synthesizing the mathematical and humanitarian cultures.

Alfred Marshall (1842--1924), the founder of the Cambridge school of neoclassicals,
``Marshallians,''
wrote in his magnum opus \cite{Marshall}:

\smallskip
\begin{itemize}
\item[]{\small\sl\indent
The function then of analysis and deduction in economics
is not to forge a few long chains of reasoning, but to forge
rightly many short chains and single connecting links\dots.
\cite[Appendix C: The Scope and Method of Economics. \S~3.]{Marshall}}
\smallskip
\item[]{\small\sl\indent
It is obvious that there is no room in economics
for long trains of deductive reasoning.
 \cite [Appendix~D: Use of Abstract Reasoning in Economics]{Marshall}}
\end{itemize}
\smallskip

\noindent
In 1906 he formulated his scepticism in regard to
mathematics as follows \cite[p.~294]{Brue}:

\smallskip
\begin{itemize}
\item[]{\small\sl\indent
[I had] a growing feeling in the
later years of my work at the subject that a good mathematical theorem
dealing with economic hypotheses was very unlikely to be good
economics: and I went more and more on the rules --- 

(1) Use mathematics
as a shorthand language, rather than an engine of inquiry. 

(2) Keep to  them till you have done. 

(3) Translate into English. 

(4) Then illustrate by examples that are important in real life. 

(5) Burn the  mathematics. 

(6) If you can't succeed in (4), burn (3). This last I
did often.} 
\end{itemize}

\noindent
Marshall intentionally counterpose the economical
and mathematical ways of thinking, noting that
the numerous short ``combs'' are appropriate in a concrete economical analysis.
Clearly, the image of a~``comb'' has nothing in common with
the upside-down pyramid, the cumulative hierarchy of
the von Neumann universe, the residence of the modern
Zermelo--Fraenkel set theory.
It is from  the times of Ancient Hellada 
that the beauty and power of mathematics rest on
the axiomatic method which presumes the derivation of new facts by 
however lengthy chains of formal implications.

The conspicuous discrepancy between economists and mathematicians
in mentality has hindered their mutual understanding and cooperation.
Many partitions, invisible  but ubiquitous, were erected in ratiocination,
isolating  the economic community from its mathematical
counterpart and vice versa.

This status quo  with deep roots in history
was always a~challenge to Kantorovich, contradicting
his views of interaction between mathematics and economics.
His path in science is well marked with  the signposts conveying
the slogan: ``Mathematicians and
Economists of the World, Unite!''
His message has been received as witnessed by the curricula and syllabi
of every economics department in a major university throughout the world.

Despite the antediluvian opinion that ``the mathematical scientist emperor of
mainstream economics is without any clothes''
(cp.~\cite{Davidson}),
the gadgets of mathematics and the idea of optimality will come in handy
for the~working economist. Calculation will supersede prophesy.
Economics as a~boon companion of mathematics will avoid merging
into any esoteric part  of~the humanities, or politics, or belles-lettres.
The new generations of mathematicians will treat the puzzling problems
of economics as an inexhaustible source of inspiration and
an attractive arena for applying and refining their formal methods.

The life of Kantorovich is a path of a scholar and  citizen whose
scientific work is inseparable from the fates of his fellows,
exemplifying the idea of serving the genuine of one's Fatherland
regardless of whatever ideological situation. This lesson
is of major import these days. Attempts at slandering or hushing up
the life and legacy of Kantorovich are doomed to demolishing.
Pygmies can never hide a giant.
}

\bibliographystyle{plain}

\begin{thebibliography}{99}


\bibitem{Kantorovich}
Kantorovich L.~V.,
{\it Selected Works. Parts 1 and 2}. Amsterdam: Gordon and Breach (1996).



\bibitem{Dantzig}
Dantzig G.~B.
{\it Linear Programming and Extensions}.
Princeton: Princeton University Press
(1963).

\bibitem{Reshetnyak}
Reshetnyak Yu.~G. and Kutateladze S.~S.``A.~Letter of N.~N.~Luzin
to L.~V.~Kantorovich,''  {\it Vestnik Ross. Acad. Nauk}, {\bf72}:8 (2002), pp.~740--742.


\bibitem{Akilov}
Kantorovich L.~V. and Akilov G.~P.
{\it Functional Analysis}. Oxford etc.: Pergamon Press,
1982.



\bibitem{IBA}
Kusraev A.~G. and Kutateladze S.~S. (2005)
{\it Boolean Valued Analysis}. Dordrecht etc.: Kluwer Academic Publishers (1999).


\bibitem{Marshall}
 Marshall A.
{\it  Principles of Economics. 8th edition}.
London: Macmillan and Co., Ltd. (1920).


\bibitem{Brue}
Brue~S.~L.
{\it The Evolution of Economic Thought. 5th edition}.
 Fort Worth: Harcourt College Publishers (1993).


\bibitem{Davidson}
Davidson P.
``Is a`Mathematical Science' an Oxymoron When Used to Describe Economics,''
{\it Post Keynesian Economics}, {\bf 25}:4 (2003).


\end{thebibliography}

\end{document}